\documentclass[12pt]{article}
\usepackage{color}
\definecolor{darkblue}{rgb}{0.00,0.25,0.50}
\usepackage[colorlinks,filecolor=blue,citecolor=darkblue]{hyperref}

\usepackage{amsfonts,amssymb,amsmath,amsthm}
\usepackage{url}
\usepackage{enumerate}
\usepackage[ukrainian, russian, english]{babel}
\usepackage[cp1251]{inputenc}
\begin{document} \selectlanguage{ukrainian}
\thispagestyle{empty}

\title{}

UDC 517.51 \vskip 5mm

\begin{center}
\textbf{\Large The best approximation of functions from anisotropic Nikol'skii--Besov classes defined in \boldmath{$\mathbb{R}^d$}}
\end{center}

\vskip 3mm

\begin{center}
\textbf{\Large  Найкраще наближення функцій з
анізотропних класів Нікольського--Бєсова визначених на \boldmath{$\mathbb{R}^d$}}
\end{center}
\vskip0.5cm

\begin{center}
 S.~Ya.~Yanchenko\\ \emph{\small
Institute of Mathematics of NAS of
Ukraine, Kyiv}
\end{center}
\begin{center}
C.~Я.~Янченко \\
\emph{\small Інститут математики НАН України, Київ}
\end{center}
\vskip0.5cm

\begin{abstract}

We establish the exact-order estimates of the best approximations of the functions from anisotropic Nikol'skii--Besov classes of several variables by entire functions in the Lebesgue spaces.
\vskip 3 mm

 Одержано точні за порядком оцінки найкращого наближення функцій з анізотропних класів
Нікольського--Бєсова функцій багатьох змінних цілими функціями у просторах Лебега.
\end{abstract}

\vskip 0.5 cm


У роботі досліджується питання найкращого наближення
функцій з анізотропних класів Нікольського--Бєсова $B^{\boldsymbol{r}}_{p,\theta}(\mathbb{R}^d)$~\cite{Nikolsky_1951}, \cite{Besov_1961}, де параметр $\boldsymbol{r}$~--- $d$-вимірний вектор з додатними координатами ($r_j>0$, $j=\overline{1,d}$). Похибка наближення при цьому вимірюється у метриці простору  $L_{q}(\mathbb{R}^d)$ $1<p\leqslant q<\infty$.  В якості апарату наближення використовуються цілі функції експоненціального типу (див., наприклад, \cite{Nikolsky_1969_book}) з носіями їх перетворення Фур'є в $d$-вимірних паралелепіпедах.

\vskip 3 mm

\textbf{1. Основні позначення та означення класів
Нікольського--Бєсова.} Нехай  $\mathbb{R}^d$, $d\geqslant 1$~---
$d$-вимірний евклідів простір з елементами
${\boldsymbol{x}=(x_1,...,x_d)}$,
${(\boldsymbol{x},\boldsymbol{y})=x_1y_1+...+x_dy_d}$.
${L_p(\mathbb{R}^d)}$, ${1\leqslant p\leqslant\infty}$,~--- простір
вимірних на $\mathbb{R}^d$ функцій ${f(\boldsymbol{x})=f(x_1,...,x_d)}$ зі скінченною
нормою
 $$
\|f\|_{L_p(\mathbb{R}^d)}=\|f\|_p:=
\left(\int\limits_{\mathbb{R}^{d}}|f(\boldsymbol{x})|^{p}d\boldsymbol{x}
\right) ^{\frac{1}{p}}, \ 1\leqslant p<\infty,
 $$
 $$
  \|f\|_{L_\infty(\mathbb{R}^d)}=\|f\|_{\infty}:=\mathop {\rm ess \sup}\limits_{\boldsymbol{x}\in \mathbb{R}^d}
  |f(\boldsymbol{x})|.
 $$

Для $k\in \mathbb{N}$, $\boldsymbol{h}\in \mathbb{R}^d$ означимо модуль гладкості $k$-го
 порядку функції ${f\in L_p(\mathbb{R}^d)}$ за змінною $x_i$,
який будемо позначати $\omega_k(f,{t}e_i)_p$, такою формулою:
$$
\omega_k(f,{t}e_i)_p= \sup \limits_{|\boldsymbol{h}|\leqslant {t}}
\|\Delta^k(f,\boldsymbol{h}e_i)\|_p=\sup \limits_{|\boldsymbol{h}|\leqslant {t}}
 \bigg\|\sum\limits_{l=0}^k(-1)^{k-l}C_k^l f(\boldsymbol{x}+l\boldsymbol{h}e_i)\bigg\|_p,
$$
де $|\boldsymbol{h}|=\sqrt{h_1^2+\ldots+h_d^2}$~--- евклідова норма вектора $\boldsymbol{h}$,
а $e_i$~--- одиничний вектор, який направлений вздовж осі $x_i$.

Нехай $r_i>0$, $r_i=\bar{r_i}+\alpha_i$, де $\bar{r_i}$~--- ціле, $0<\alpha_i\leqslant 1$, $i=\overline{1,d}$.

\vskip 2 mm \textbf{\emph{Означення} 1.} Будемо говорити, що функція $f \in L_p(\mathbb{R}^d)$ належить
простору $B^{\boldsymbol{r}}_{p,\theta}(\mathbb{R}^d)$,
$1\leqslant p, \ \theta\leqslant\infty$, $\boldsymbol{r}>0$, якщо вона має інтегровані  в степені $p$  на
$\mathbb{R}^d$  часткові, узагальнені в сенсі Соболєва, похідні вигляду
$$
D_i^k f=\frac{\partial^k f}{\partial^k x_i}, \ k=\overline{0,\bar{r_i}}, \ i=\overline{1,d},
$$
і при цьому
$$
\|f\|_{B^{\boldsymbol{r}}_{p,\theta}}=\|f\|_p+\sum\limits_{i=1}^d\left(\int\limits^{\infty}_0
{t}^{-\theta \alpha_i-1}
\omega_{1+[\alpha_i]}^{\theta}(D_i^{\bar{r_i}}f,{t}e_i)_p {dt}\right)^{\frac{1}{\theta}} < \infty \ \mbox {при} \  1 \leqslant  \theta < \infty
$$
та
$$
\|f\|_{B^{\boldsymbol{r}}_{p,\infty}}=\|f\|_{H^{\boldsymbol{r}}_{p}}=\|f\|_{\infty}+\sum\limits_{i=1}^d \sup
{t}^{-\alpha_i-1}
\omega_{1+[\alpha_i]}(D_i^{\bar{r_i}}f,{t}e_i)_p  < \infty \ \mbox
{при} \ \theta = \infty.
$$

Відзначимо, що з так введеною нормою простори $B^{\boldsymbol{r}}_{p,\theta}(\mathbb{R}^d)$ будуть банаховими.

Простори $B^{\boldsymbol{r}}_{p,\theta}(\mathbb{R}^d)$ були введені О.\,В.~Бєсовим~\cite{Besov_1961} і
$B^{\boldsymbol{r}}_{p,\infty}(\mathbb{R}^d)=H^{\boldsymbol{r}}_p(\mathbb{R}^d)$, де
$H^{\boldsymbol{r}}_p(\mathbb{R}^d)$~--- простори, які ввів
 С.\,М.~Нікольський~\cite{Nikolsky_1951}. Далі, зберігаючи ті самі позначення, будемо розглядати класи $B^{\boldsymbol{r}}_{p,\theta}(\mathbb{R}^d)$, тобто одиничні кулі у просторах $B^{\boldsymbol{r}}_{p,\theta}(\mathbb{R}^d)$:
$$
B^{\boldsymbol{r}}_{p,\theta}(\mathbb{R}^d):= \left\{f\in L_p(\mathbb{R}^d):
\|f\|_{B^{\boldsymbol{r}}_{p,\theta}(\mathbb{R}^d)}\leqslant 1\right\}.
$$
Крім цього, для спрощення записів, замість $B^{\boldsymbol{r}}_{p,\theta}(\mathbb{R}^d)$ та
$H^{\boldsymbol{r}}_{p}(\mathbb{R}^d)$ будемо використовувати  позначення
$B^{\boldsymbol{r}}_{p,\theta}$ та $H^{\boldsymbol{r}}_{p}$.

Зазначимо, що важливим для встановлення результатів є той
факт, що простори $B^{\boldsymbol{r}}_{p,\theta}$ зі зростанням параметра $\theta$
розширюються (див., наприклад, \cite[с.~278]{Nikolsky_1969_book}),
тобто
\begin{equation}\label{vklad}
 B^{\boldsymbol{r}}_{p,1}\subset B^{\boldsymbol{r}}_{p,\theta}\subset B^{\boldsymbol{r}}_{p,\theta'}\subset
 B^{\boldsymbol{r}}_{p,\infty}=H^{\boldsymbol{r}}_p, \ \ \ 1\leqslant\theta<\theta' \leqslant
 \infty.
\end{equation}

Наведемо один результат П.\,І.~Лізоркіна~(див. \cite{Lizorkin_1968_sib}), який
дає можливість означити норму функцій з просторів $B^{\boldsymbol{r}}_{p,\theta}(\mathbb{R}^d)$ в
іншій формі, яка в подальшому зумовлює використання перетворення Фур'є в теорії даних просторів.

Для цього попередньо наведемо необхідні означення.

Назвемо найкращим  наближенням функції $f\in L_p(\mathbb{R}^d)$ за допомогою цілих функцій
степеня $\nu_1,\ldots,\nu_d$ величину
\begin{equation}\label{En_best}
E_{\nu_1,\ldots,\nu_d}(f)_p=\inf\limits_{g_{\nu_1,\ldots,\nu_d}}\|f-g_{\nu_1,\ldots,\nu_d}\|_p
\end{equation}
де $\inf$ береться по всіх цілих функціях $g_{\nu_1,\ldots,\nu_d}(x_1,\ldots,x_d)$
степеня $\nu_1,\ldots,\nu_d$ відповідно за змінними $x_1,\ldots,x_d$.

\vskip 2 mm

\bf Теорема A \rm\cite{Lizorkin_1968_sib}\textbf{.}  \it  Функція $f$ належить простору $B^{\boldsymbol{r}}_{p,\theta}(\mathbb{R}^d)$, $\boldsymbol{r}>0$,
${1\leqslant p,\theta\leqslant \infty }$, тоді і тільки тоді, коли вона зображується збіжним у метриці простору $L_p(\mathbb{R}^d)$ рядом
\begin{equation}\label{f-efet}
 f(\boldsymbol{x})=\sum\limits_{s=0}^{\infty}P_{\boldsymbol{a}^s}(\boldsymbol{x}),
 \ \ P_{\boldsymbol{a}^s}(\boldsymbol{x})=P_{a_1^s,\ldots,a_d^s}(\boldsymbol{x}),
\end{equation}
де $P_{\nu_1,\ldots,\nu_d}(\boldsymbol{x})$~--- цілі функції степеня не вищого за
 $\nu_1,\ldots,\nu_d$ по кожній змінній $x_1,\ldots, x_d$  відповідно, і виконується умова
\begin{equation}\label{f-norm-dek1}
 \left(\sum\limits_{s=0}^{\infty}b^{s\theta}\|P_{\boldsymbol{a}^s}\|_p^{\theta}
 \right)^{\frac{1}{\theta}}<\infty,  \ \ \mbox{де} \ \ b=a_i^{r_i}>1, \ i=\overline{1,d}.
\end{equation}
Окрім цього має місце оцінка
\begin{equation}\label{f-norm-dek2}
 \|f\|_{B^{\boldsymbol{r}}_{p,\theta}}\leqslant C_1\left(\sum\limits_{s=0}^{\infty}b^{s\theta}\|P_{\boldsymbol{a}^s}\|_p^{\theta}
 \right)^{\frac{1}{\theta}}.
\end{equation}

Якщо, крім того, частинні суми $n$-го порядку  ряду $(\ref{f-efet})$ реалізують
найкраще наближення або дають порядок найкращого наближення, то вираз у лівій частині
 $(\ref{f-norm-dek1})$ і $\|f\|_{B^{\boldsymbol{r}}_{p,\theta}}$ еквівалентні,
  тобто разом із $\eqref{f-norm-dek2}$ має місце оцінка
$$
\left(\sum\limits_{s=0}^{\infty}b^{s\theta}\|P_{\boldsymbol{a}^s}\|_p^{\theta}
 \right)^{\frac{1}{\theta}}\leqslant C_2\|f\|_{B^{\boldsymbol{r}}_{p,\theta}} .
$$
\rm \vskip 2 mm

На основі  теореми~А дамо еквівалентне означення анізотропних просторів $B^{\boldsymbol{r}}_{p,\theta}$, яким будемо
користуватися у подальших міркуваннях. Для цього
нагадаємо означення перетворення Фур'є (див., наприклад,
\cite{Lizorkin_69})  з використанням якого дається
відповідне означення.

Нехай $S=S(\mathbb{R}^d)$~--- простір Л.~Шварца основних нескінченно
диференційовних на $\mathbb{R}^d$ комплекснозначних функцій
$\varphi$, що спадають на нескінченності разом зі своїми похідними
швидше за будь-який степінь функції $|\boldsymbol{x}|^{-1}$ (див., наприклад,
\cite{Lizorkin_69}, \cite{Vladimirov}). Через $S'$ позначимо
простір лінійних неперервних функціоналів на $S$. Зазначимо, що
елементами простору $S'$ є узагальнені функції. Якщо $f\in S'$ і
$\varphi \in S$, то $\langle f,\varphi\rangle$ позначає значення $f$
на $\varphi$.

Перетворення Фур'є $\mathfrak{F}\varphi: S\rightarrow S$
визначається згідно з формулою
$$
(\mathfrak{F}\varphi)(\boldsymbol{\lambda})=\frac{1}{(2\pi)^{d/2}}\int
 \limits_{\mathbb{R}^d}\varphi(\boldsymbol{t})
 e^{-i(\boldsymbol{\lambda},\boldsymbol{t})}d\boldsymbol{t}
 \equiv \widetilde{\varphi}(\boldsymbol{\lambda}).
$$

Обернене перетворення Фур'є $\mathfrak{F}^{-1}\varphi:\ S\rightarrow S$ задається таким чином:
$$
(\mathfrak{F}^{-1}\varphi)(\boldsymbol{t})=\frac{1}{(2\pi)^{d/2}}
\int \limits_{\mathbb{R}^d}\varphi(\boldsymbol{\lambda})
e^{i(\boldsymbol{\lambda},\boldsymbol{t})}d\boldsymbol{\lambda}\equiv
\widehat{\varphi}(\boldsymbol{t}).
$$

Перетворення Фур'є узагальнених функцій $f\in S'$ (для нього ми
зберігаємо те ж позначення)  визначається згідно з формулою
$$
\langle \mathfrak{F}f,\varphi\rangle=\langle f,\mathfrak{F}\varphi
\rangle  \ \ \ (\langle \widetilde{f},\varphi\rangle=\langle
f,\widetilde{\varphi} \rangle),
$$
де $\varphi \in S$.

Обернене перетворення Фур'є узагальненої функції  $f\in S'$ також позначимо
$\mathfrak{F}^{-1}f$,  і визначається воно аналогічно до прямого
перетворення Фур'є згідно з формулою
$$
\langle \mathfrak{F}^{-1}f,\varphi\rangle=\langle
f,\mathfrak{F}^{-1}\varphi \rangle  \ \ \ (\langle
\widehat{f},\varphi\rangle=\langle f,\widehat{\varphi} \rangle).
$$

Зазначимо, що кожна функція $f \in L_p(\mathbb{R}^d)$, $1\leqslant p \leqslant \infty$,
 визначає лінійний неперервний функціонал на $S$ згідно з формулою
$$
\langle f,\varphi \rangle = \int \limits_{\mathbb{R}^d}
 f(\boldsymbol{x})\varphi(\boldsymbol{x}) d\boldsymbol{x}, \ \ \varphi\in S,
$$
і, як наслідок, у цьому сенсі вона є елементом $S'$. Тому
перетворення Фур'є функції $f \in L_p(\mathbb{R}^d)$, $1\leqslant p \leqslant \infty$, можна
розглядати як перетворення Фур'є узагальненої функції $\langle f,\varphi \rangle$.

Носієм узагальненої функції $f$ будемо називати замикання
$\overline{\mathfrak{N}}$ такої множини точок
$\mathfrak{N}\subset\mathbb{R}^d$, що для довільної $\varphi \in S$,
яка дорівнює нулю в $\overline{\mathfrak{N}}$, виконується рівність
$\langle f,\varphi \rangle = 0$. Носій узагальненої функції $f$
будемо позначати через $\mbox{supp}\, f$. Також будемо говорити, що функція $f$
 зосереджена на множині $G$, якщо $\mbox{supp}\, f \subseteq G$.

У подальшому будемо користуватися такими позначеннями. Нехай функція $f$ представлена інтегралом Фур'є
$$
f(\boldsymbol{x})=\frac{1}{(2\pi)^{d/2}}
\int\limits_{\mathbb{R}^d}\tilde{f}(\boldsymbol{\lambda})e^{i(\boldsymbol{\lambda}, \boldsymbol{x})}d\boldsymbol{\lambda}.
$$
Тоді відрізком інтегралу Фур'є функції $f$  назвемо вираз
$$
S_{\boldsymbol{\sigma}}(f)=\frac{1}{(2\pi)^{d/2}}\int\limits_{-\sigma_1}^{\sigma_1}\ldots
\int\limits_{-\sigma_d}^{\sigma_d}
\tilde{f}(\boldsymbol{\lambda})e^{i(\boldsymbol{\lambda}, \boldsymbol{x})}d\boldsymbol{\lambda},
$$
де $\tilde{f}(\boldsymbol{\lambda})$~--- перетворення Фур'є функції $f\in L_p(\mathbb{R}^d)$.

Нехай $D_{\boldsymbol{a}^s}=D_{a_1^{s},\ldots, a_d^{s}}$~--- паралелепіпед:
$|\lambda_j|<a_j^s$, $j=\overline{1,d}$, $s\geqslant 0$, а ${\Gamma_{\boldsymbol{a}^s}=D_{\boldsymbol{a}^s}-D_{\boldsymbol{a}^{s-1}}}$ при $s\geqslant 1$ і $\Gamma_{\boldsymbol{a}^0}=D_{\boldsymbol{a}^0}$. Покладемо
$$
f_{\boldsymbol{a}^s}=S_{\boldsymbol{a}^s}(f)-S_{\boldsymbol{a}^{s-1}}(f)= \int\limits_{\Gamma_{\boldsymbol{a}^s}}\tilde{f}(\boldsymbol{\lambda})e^{i(\boldsymbol{\lambda}, \boldsymbol{x})}d\boldsymbol{\lambda}, s\geqslant 1,
$$
і
$$
f_{\boldsymbol{a}^0}=S_{\boldsymbol{a}^0}(f)= \int\limits_{\Gamma_{\boldsymbol{a}^0}}\tilde{f}(\boldsymbol{\lambda})e^{i(\boldsymbol{\lambda}, \boldsymbol{x})}d\boldsymbol{\lambda}.
$$
Представлення функції $f$ рядом
$$
f=f_{\boldsymbol{a}^0}+\sum\limits_{s=1}^{\infty}f_{\boldsymbol{a}^s}= \sum\limits_{s=0}^{\infty}f_{\boldsymbol{a}^s}
$$
будемо називати розшаруванням $f$  ($\boldsymbol{a}$-розшаруванням $f$). У випадку, коли $f\in L_p$, ${p>2}$, $S_{\boldsymbol{a}^s}(f)$ розуміють, взагалі кажучи, як результат дії на $f$ деякого оператора, який в образах Фур'є зводиться до множення на характеристичну функцію області $D_{\boldsymbol{a}^s}$ (див. \cite{Lizorkin_69}, (\S3, гл.1)).

Для функції $f\in L_p(\mathbb{R}^d)$ розглянемо величину
\begin{equation}\label{En_Fourier}
\mathcal{E}_{D_{\boldsymbol{a}^n}}(f)_p=\|f-S_{\boldsymbol{a}^{n-1}}(f)\|_p, \ \ n\in \mathbb{N},
\end{equation}
яка називається наближенням функції $f$  $\boldsymbol{a}^n$-відрізками інтеграла Фур'є.

У випадку $1<p<\infty$ величини (\ref{En_best}) і (\ref{En_Fourier}) мають один і той же порядок~\cite{Lizorkin_1968_sib}, тобто для функції $f\in L_p(\mathbb{R}^d)$
\begin{equation}\label{En_best-En_Fourier}
\mathcal{E}_{D_{\boldsymbol{a}^n}}(f)_p\asymp E_{\boldsymbol{a}^n}(f)_p.
\end{equation}

Далі для вектора $\boldsymbol{r}=(r_1,\ldots,r_d)$, $r_j>0$, $j=\overline{1,d}$, введемо величину
\begin{equation}\label{g(r)}
g(\boldsymbol{r})=\Bigg(\frac{1}{d}\sum\limits_{j=1}^{d}\frac{1}{r_j}\Bigg)^{-1}.
\end{equation}

Зауважимо, що при $r_1=r_2=\ldots=r_d=r$ маємо $g(\boldsymbol{r})=r$.

Тоді для норми анізотропних просторів $B^{\boldsymbol{r}}_{p,\theta}(\mathbb{R}^d)$, згідно з теоремою А, можна записати співвідношення~\cite{Lizorkin_1968_sib}:
\begin{equation}\label{f-norm-dek}
 \|f\|_{B^{\boldsymbol{r}}_{p,\theta}(\mathbb{R}^d)}\asymp
 \left(\sum\limits_{s=0}^{\infty}b^{s\theta}\|f_{\boldsymbol{a}^s}\|_p^{\theta}
 \right)^{\frac{1}{\theta}}<\infty,   \ \ \mbox{при} \ \ 1\leqslant \theta <\infty,
\end{equation}
\begin{equation}\label{f-norm-dek-inf}
 \|f\|_{B^{\boldsymbol{r}}_{p,\infty}(\mathbb{R}^d)}\asymp
 \sup\limits_{s\geqslant 0} b^{s}\|f_{\boldsymbol{a}^s}\|_p <\infty,
\end{equation}
де $b=2^{g(\boldsymbol{r})}$, тобто $a_j=2^{g(\boldsymbol{r})/r_j}$, $j=\overline{1,d}$.

\vskip 5 mm

\textbf{2. Наближення класів \boldmath{$B^{\boldsymbol{r}}_{p,\theta}(\mathbb{R}^d)$}  у метриці простору \boldmath{$L_q(\mathbb{R}^d)$, ${1<p\leqslant q< \infty}$}.}

Попередньо сформулюємо твердження, яке буде істотно використовуватися при
встановленні результатів.

\vskip 1 mm
\bf Теорема Б \rm \cite[c.~150]{Nikolsky_1969_book}. \it  Якщо
$1\leqslant p_1 \leqslant p_2 \leqslant \infty$, то для цілої функції
експоненціального типу $g=g_{\boldsymbol{\nu}}\in L_p(\mathbb{R}^d)$ має місце
``нерівність різних метрик''
\begin{equation}\label{Riz_Metric}
 \|g_{\boldsymbol{\nu}}\|_{L_{p_2}(\mathbb{R}^d)}\leqslant 2^d\left( \prod \limits_{j=1}^d
 \nu_k\right)^{\frac{1}{p_1}-\frac{1}{p_2}}\|g_{\boldsymbol{\nu}}\|_{L_{p_1}(\mathbb{R}^d)}.
\end{equation}\rm \vskip 1 mm

Наведемо одержані результати.

\vskip 1 mm \bf Теорема 1. \it Нехай $1<p \leqslant q<\infty$, $1\leqslant \theta
\leqslant \infty$. Тоді для  ${g(\boldsymbol{r})>d\left(\frac{1}{p}-\frac{1}{q}\right)}$, мають місце порядкові співвідношення
\begin{equation} \label{teor_E_n_aniz}
   {\mathcal{E}}_{D_{\boldsymbol{a}^n}}(B^{\boldsymbol{r}}_{p,\theta})_{q}\asymp {{E}}_{\boldsymbol{a}^n}(B^{\boldsymbol{r}}_{p,\theta})_{q}
   \asymp 2^{-n\left(g(\boldsymbol{r})-d\left(\frac{1}{p}-\frac{1}{q}\right)\right)},
\end{equation}
де $a_j=2^{g(\boldsymbol{r})/r_j}$, $j=\overline{1,d}$.

{\textbf{\textit{Доведення.}}} \ \rm Спочатку отримаємо
 в (\ref{teor_E_n_aniz}) оцінки зверху. Оскільки ${B^{\boldsymbol{r}}_{p,\theta}\subset
 B^{\boldsymbol{r}}_{p,\infty}=H^{\boldsymbol{r}}_p}$, $1\leqslant \theta <\infty$, то шукану оцінку достатньо
 отримати для величини $\mathcal{E}_{D_{\boldsymbol{a}^n}}(H_p^{\boldsymbol{r}})_q$. У залежності від
 співвідношення між параметрами $p$ i $q$  розглянемо два
 випадки.

 1) Нехай  $1<p=q<\infty$. Оскільки для $f\in H^{\boldsymbol{r}}_p$ згідно з  (\ref{f-norm-dek-inf})  ${\|f_{\boldsymbol{a}^s}\|_p\ll 2^{-sg(\boldsymbol{r})}}$, то  скориставшись нерівністю Мінковського, будемо мати
$$
{\mathcal{E}}_{D_{\boldsymbol{a}^n}}(f)_{p}=\|f-S_{\boldsymbol{a}^{n-1}}(f)\|_p = \Big\|\sum
 \limits_{s=0}^{\infty}f_{\boldsymbol{a}^s}-S_{\boldsymbol{a}^{n-1}}(f)\Big\|_p=
$$
 $$
   = \Big\|\sum
 \limits_{s=n}^{\infty}f_{\boldsymbol{a}^s}\Big\|_p\leqslant \sum
 \limits_{s=n}^{\infty}\|f_{\boldsymbol{a}^s}\|_p \leqslant \sum
 \limits_{s=n}^{\infty} 2^{-sg(\boldsymbol{r})}\ll 2^{-ng(\boldsymbol{r})}.
 $$

 2) Нехай тепер $1<p<q<\infty$. Тоді для $f\in H^{\boldsymbol{r}}_p$, врахувавши (\ref{f-norm-dek-inf}) та скориставшись
 нерівностями Мінковського і  різних метрик
 Нікольського (\ref{Riz_Metric}), можемо записати
 $$
 \mathcal{E}_{D_{\boldsymbol{a}^n}}(f)_{q}=\|f-S_{\boldsymbol{a}^{n-1}}(f)\|_q=\Big\|\sum
 \limits_{s=n}^{\infty}f_{\boldsymbol{a}^s}\Big\|_q\leqslant \sum
 \limits_{s=n}^{\infty}\|f_{\boldsymbol{a}^s}\|_q\ll
 $$
 $$
\ll \sum
 \limits_{s=n}^{\infty}2^{sd\left(\frac{1}{p}-\frac{1}{q}\right)}\|f_{\boldsymbol{a}^s}\|_p\ll \sum
 \limits_{s=n}^{\infty}2^{sd\left(\frac{1}{p}-\frac{1}{q}\right)} 2^{-sg(\boldsymbol{r})}=
 $$
 $$
  = \sum
 \limits_{s=n}^{\infty}2^{-s\left(g(\boldsymbol{r})-d\left(\frac{1}{p}+\frac{1}{q}\right)\right)}\ll
 2^{-n\left(g(\boldsymbol{r})-d\left(\frac{1}{p}-\frac{1}{q}\right)\right)}.
 $$

Оцінку зверху для величини $\mathcal{E}_{D_{\boldsymbol{a}^n}}(H_p^{\boldsymbol{r}})_{q}$ і, таким чином, згідно з (\ref{En_best-En_Fourier}) для найкращого наближення $E_{D_{\boldsymbol{a}^n}}(H_p^{\boldsymbol{r}})_{q}$ встановлено.

Отримаємо тепер в (\ref{teor_E_n_aniz}) оцінки знизу. Оскільки має
місце вкладення ${B^{\boldsymbol{r}}_{p,1}\subset B^{\boldsymbol{r}}_{p,\theta}}$,
$1<\theta\leqslant\infty$, то шукану оцінку  достатньо отримати для величини
$\mathcal{E}_{D_{\boldsymbol{a}^n}}(B^{\boldsymbol{r}}_{p,1})_q$. Іншими словами достатньо оцінити знизу
величину ${\|f-S_{\boldsymbol{a}^{n-1}}(f)\|_q}$ для деякої функції ${f\in
B^{\boldsymbol{r}}_{p,1}}$.

З цією метою розглянемо функцію  $F_{\boldsymbol{k}}(\boldsymbol{x})$ на основі якої побудуємо функцію для якої досягається оцінка (\ref{teor_E_n_aniz}).

Нехай $\boldsymbol{k}\in \mathbb{N}^d$, $\boldsymbol{k}=(k,\ldots,k)$,
$$
 F_{\boldsymbol{k}}(\boldsymbol{x})=\prod\limits_{j=1}^{d}\sqrt{\frac {2}{\pi}}\ \frac{ \sin{a_j^k x_j}}{x_j}- \prod\limits_{j=1}^{d}\sqrt{\frac {2}{\pi}}\ \frac{ \sin{a_j^{k-1} x_j}}{x_j}
 $$
 та
 $$
F_{0}(\boldsymbol{x})=\prod\limits_{j=1}^{d}\sqrt{\frac {2}{\pi}}\ \frac{ \sin{x_j}}{x_j}.
 $$

 Тоді для перетворення Фур'є функції $F_{\boldsymbol{k}}(\boldsymbol{x})$
 має місце співвідношення
 $$
 \mathfrak{F}F_{\boldsymbol{k}}(\boldsymbol{x})=\chi_{\boldsymbol{k}}(\boldsymbol{\lambda})=\prod \limits_{j=1}^d \chi_{k}(\lambda_j),
 $$
де

 \begin{minipage}{9 cm}
    $$
      \chi_{k}(\lambda_j)=
 \begin{cases}
    1, & a_j^{k-1}<|\lambda_j|<a_j^{k}, \\
    \frac{1}{2}, & |\lambda_j|=a_j^{k-1} \ \mbox{або} \ |\lambda_j|=a_j^{k}, \\
    0 & \mbox{--- в інших випадках},
 \end{cases}
    $$
\end{minipage}
\begin{minipage}{6 cm}
  $$
   \chi_{0}(x_j)=
 \begin{cases}
    1, & |\lambda_j|<1; \\
    \frac{1}{2}, & |\lambda_j|=1; \\
    0, & |\lambda_j|>1.
 \end{cases}
  $$
\end{minipage}
\vskip 1mm

Відповідно для оберненого перетворення будемо мати
 $$
 \mathfrak{F}^{-1}\chi_{\boldsymbol{k}}(\boldsymbol{\lambda})=F_{\boldsymbol{k}}(\boldsymbol{x}).
 $$

Відзначимо, що таким чином  $F_{\boldsymbol{k}}(\boldsymbol{x})$~--- ціла функція з $L_p(\mathbb{R}^d)$, носій перетворення Фур'є якої зосереджене в $\Gamma_{\boldsymbol{a}^{\boldsymbol{k}}}$.

Перш ніж безпосередньо перейти до встановлення оцінки знизу в (\ref{teor_E_n_aniz}),  одержимо порядок величини
\begin{equation}\label{F_k_norm}
\|F_{\boldsymbol{k}}\|_p=\bigg\|\prod\limits_{j=1}^{d}\sqrt{\frac {2}{\pi}}\ \frac{ \sin{a_j^k x_j}}{x_j}- \prod\limits_{j=1}^{d}\sqrt{\frac {2}{\pi}}\ \frac{ \sin{a_j^{k-1} x_j}}{x_j}\bigg\|_p.
\end{equation}

Для оцінки зверху будемо мати
$$
\|F_{\boldsymbol{k}}\|_p= \bigg\|\prod\limits_{j=1}^{d}\sqrt{\frac {2}{\pi}}\ \frac{ \sin{a_j^k x_j}}{x_j} - \prod\limits_{j=1}^{d}\sqrt{\frac {2}{\pi}}\ \frac{ \sin{a_j^{k-1} x_j}}{x_j} \bigg\|_p\leqslant
$$
$$
\leqslant \bigg\|\prod\limits_{j=1}^{d}\sqrt{\frac {2}{\pi}}\ \frac{ \sin{a_j^k x_j}}{x_j}\bigg\|_p + \bigg\|\prod\limits_{j=1}^{d}\sqrt{\frac {2}{\pi}}\ \frac{ \sin{a_j^{k-1} x_j}}{x_j} \bigg\|_p =
$$
$$
 = \left(\int\limits_{\mathbb{R}^d} \Big|\prod\limits_{j=1}^{d}\sqrt{\frac {2}{\pi}}\ \frac{ \sin{a_j^k x_j}}{x_j}\Big|^p \prod\limits_{j=1}^{d}dx_j\right)^{\frac{1}{p}} + \left(\int\limits_{\mathbb{R}^d} \Big| \prod\limits_{j=1}^{d}\sqrt{\frac {2}{\pi}}\ \frac{ \sin{a_j^{k-1} x_j}}{x_j}\Big|^p \prod\limits_{j=1}^{d} dx_j\right)^{\frac{1}{p}}=
$$
$$
 = \left(\frac {2}{\pi}\right)^{\frac{d}{2}}\prod\limits_{j=1}^{d}\left(\int\limits_{\mathbb{R}^d} \Big| \frac{ \sin{a_j^k x_j}}{x_j}\Big|^p dx_j\right)^{\frac{1}{p}} + \left(\frac {2}{\pi}\right)^{\frac{d}{2}}\prod\limits_{j=1}^{d}\left(\int\limits_{\mathbb{R}^d} \Big| \frac{ \sin{a_j^{k-1} x_j}}{x_j}\Big|^p dx_j\right)^{\frac{1}{p}}=
$$
$$
=\left(\frac {2}{\pi}\right)^{\frac{d}{2}}\prod\limits_{j=1}^{d}\left( a_j^{k(p-1)} \int\limits_{\mathbb{R}^d} \Big| \frac{ \sin{x_j}}{x_j}\Big|^p dx_j\right)^{\frac{1}{p}} +
\left(\frac {2}{\pi}\right)^{\frac{d}{2}}\prod\limits_{j=1}^{d}\left( a_j^{(k-1)(p-1)} \int\limits_{\mathbb{R}^d} \Big| \frac{ \sin{x_j}}{x_j}\Big|^p dx_j\right)^{\frac{1}{p}}\ll
$$
\begin{equation}\label{F_k_norm2}
\ll \prod\limits_{j=1}^{d} a_j^{\frac{k(p-1)}{p}} + \prod\limits_{j=1}^{d} a_j^{\frac{(k-1)(p-1)}{p}}\ll \prod\limits_{j=1}^{d} a_j^{k\left(1-\frac{1}{p}\right)}=\prod\limits_{j=1}^{d} a_j^{k/p'}.
\end{equation}

Врахувавши, що $a_j=2^{g(\boldsymbol{r})/r_j}$ та співвідношення (\ref{g(r)}), оцінку (\ref{F_k_norm2}) продовжимо таким чином
\begin{equation}\label{F_k_norm_sup}
\|F_{\boldsymbol{k}}\|_p \ll \prod\limits_{j=1}^{d} a_j^{k/p'}=\prod\limits_{j=1}^{d} 2^{kg(\boldsymbol{r})/r_jp'}= 2^{\frac{kg(\boldsymbol{r})}{p'}\sum\limits_{j=1}^{d}\frac{1}{r_j}}=
2^{\frac{kg(\boldsymbol{r})}{p'}\frac{d}{g(\boldsymbol{r})}}=2^{\frac{dk}{p'}}.
\end{equation}

При оцінці норми $\|F_{\boldsymbol{k}}\|_p$ знизу, отримаємо
$$
\|F_{\boldsymbol{k}}\|_p = \bigg\|\prod\limits_{j=1}^{d}\sqrt{\frac {2}{\pi}}\ \frac{ \sin{a_j^k x_j}}{x_j} - \prod\limits_{j=1}^{d}\sqrt{\frac {2}{\pi}}\ \frac{ \sin{a_j^{k-1} x_j}}{x_j} \bigg\|_p\geqslant
$$
$$
\geqslant \bigg| \Big\|\prod\limits_{j=1}^{d}\sqrt{\frac {2}{\pi}}\ \frac{ \sin{a_j^k x_j}}{x_j}\Big\|_p - \Big\|\prod\limits_{j=1}^{d}\sqrt{\frac {2}{\pi}}\ \frac{ \sin{a_j^{k-1} x_j}}{x_j} \Big\|_p \bigg| =
$$
$$
 = \left|\left(\int\limits_{\mathbb{R}^d} \Big|\prod\limits_{j=1}^{d}\sqrt{\frac {2}{\pi}}\ \frac{ \sin{a_j^k x_j}}{x_j}\Big|^p \prod\limits_{j=1}^{d}dx_j\right)^{\frac{1}{p}} - \left(\int\limits_{\mathbb{R}^d} \Big| \prod\limits_{j=1}^{d}\sqrt{\frac {2}{\pi}}\ \frac{ \sin{a_j^{k-1} x_j}}{x_j}\Big|^p \prod\limits_{j=1}^{d} dx_j\right)^{\frac{1}{p}}\right|=
$$
$$
 = \left(\frac {2}{\pi}\right)^{\frac{d}{2}}\left|\prod\limits_{j=1}^{d}\left(\int\limits_{\mathbb{R}^d} \Big| \frac{ \sin{a_j^k x_j}}{x_j}\Big|^p dx_j\right)^{\frac{1}{p}} - \prod\limits_{j=1}^{d}\left(\int\limits_{\mathbb{R}^d} \Big| \frac{ \sin{a_j^{k-1} x_j}}{x_j}\Big|^p dx_j\right)^{\frac{1}{p}}\right|\gg
$$
\begin{equation}\label{F_k_norm_inf}
\gg \left|\prod\limits_{j=1}^{d} a_j^{\frac{k(p-1)}{p}} - \prod\limits_{j=1}^{d} a_j^{\frac{(k-1)(p-1)}{p}}\right| \gg (2^{\frac{dk}{p'}}- 2^{\frac{d(k-1)}{p'}})\gg 2^{\frac{dk}{p'}}
\end{equation}

Співставивши (\ref{F_k_norm_sup}) і (\ref{F_k_norm_inf}), для $\|F_{\boldsymbol{k}}\|_p$ можемо записати порядкове співвідношення
\begin{equation}\label{F_k_norm_os}
\|F_{\boldsymbol{k}}\|_p\asymp 2^{\frac{dk}{p'}}.
\end{equation}

Далі розглянемо функцію
$$
g_1(\boldsymbol{x})=C_1 2^{-n\left(g(\boldsymbol{r})+\frac{d}{p'}\right)}F_{\boldsymbol{n}}(\boldsymbol{x}),
$$
де $\boldsymbol{n}=(n,\ldots,n)\in \mathbb{N}^d$, $1/p+ 1/p'=1$, $C_1>0$.

Покажемо, що з деякою константою $C_1>0$ функція $g_1$ належить класу $B^{\boldsymbol{r}}_{p,1}(\mathbb{R}^d)$. Згідно з (\ref{f-norm-dek}) та (\ref{F_k_norm_os}) маємо
$$
\|g_1\|_{B^{\boldsymbol{r}}_{p,1}}\asymp \sum\limits_{s} 2^{sg(\boldsymbol{r})} \|f_{\boldsymbol{a}^s}(g_1)\|_p\asymp
$$
$$
\asymp  \sum\limits_{s} 2^{sg(\boldsymbol{r})} 2^{-n\left(g(\boldsymbol{r})+\frac{d}{p'}\right)}\|F_{\boldsymbol{n}}\|_p= 2^{-n\left(g(\boldsymbol{r})+\frac{d}{p'}\right)} \sum\limits_{s} 2^{sg(\boldsymbol{r})} 2^{\frac{dn}{p'}}\ll
$$
$$
\ll2^{-n\left(g(\boldsymbol{r})+\frac{d}{p'}\right)} 2^{ng(\boldsymbol{r})} 2^{\frac{dn}{p'}} =1.
$$

Оскільки, за вибором функції $g_1$ для неї має місце співвідношення $S_{\boldsymbol{a}^{n-1}}(g_1)=0$, то
$$
\mathcal{E}_{D_{\boldsymbol{a}^{n}}}(B^{\boldsymbol{r}}_{p,1})_q\geqslant
\mathcal{E}_{D_{\boldsymbol{a}^{n}}}(g_1)_q=\|g_1-S_{\boldsymbol{a}^{n-1}}(g_1)\|_q= \|g_1\|_q\gg
$$
$$
\gg2^{-n\left(g(\boldsymbol{r})+\frac{d}{p'}\right)}\|F_{\boldsymbol{n}}\|_q\gg 2^{-n\left(g(\boldsymbol{r})+\frac{d}{p'}\right)} 2^{\frac{dn}{q'}}=2^{-n\left(g(\boldsymbol{r})-d\left(\frac{1}{p}-\frac{1}{q}\right)\right)}.
$$

Оцінку знизу в  (\ref{teor_E_n_aniz}) встановлено. Теорему доведено.

\vskip 2 mm {\textbf{\emph{Зауваження}~1.}}  У випадку $r_1=\ldots=r_d=r$, тобто для ізотропних класів Нікольського--Бєсова $B^r_{p,\theta}(\mathbb{R}^d)$, оцінку (\ref{teor_E_n_aniz}) встановлено в \cite{Yanchenko_Zb_2010TF}. Зазначимо також, що апроксимативні характеристики ізотропних класів $B^r_{p,\theta}(\mathbb{R}^d)$ досліджувалися  у роботі~\cite{Yanchenko_2015UMG}.

\vskip 2 mm {\textbf{\emph{Зауваження}~2.}}  Анізотропні класи Нікольського--Бєсова функцій багатьох змінних, що визначені на $\mathbb{R}^d$ з точки зору знаходження точних за порядком значень деяких апроксимативних характеристик досліджувалися, зокрема, у роботах~\cite{Jiang_Yanjie_Liu_Yongping_JAT_2000}, \cite{Jiang_Yanjie_AM_2002}, а  ізотропні та анізотропні класи Нікольського--Бєсова періодичних функцій багатьох змінних~---   у роботах \cite{Romanyuk_UMG_2009}\,--\,\cite{Myronyuk_UMG_2014_8}.

\vskip 2 mm {\textbf{\emph{Зауваження}~3.}}  В одномірному  випадку $(d=1)$ анізотропні класи  Нікольського--Бєсова збігаються з класами Нікольського--Бєсова мішаної гладкості, які досліджувалися в роботах  \cite{WangHeping_SunYongsheng_1995}, \cite{Yanchenko_YMG_2010_8}. Розв'язку ряду екстремальних проблем апроксимації  функцій визначених на прямій  присвячені роботи С.\,Б.~Вакарчука~\cite{Vakarchuk_2012_1}, \cite{Vakarchuk_2012_2}, де також проведено детальний порівняльний аналіз завершених результатів, які пов'язані з розв'язком екстремальних задач теорії наближення в періодичному випадку і випадку всієї дійсної осі.

\vskip 3 mm

\textbf{Contact information:}
Department of the Theory of Functions, Institute of Mathematics of National
Academy of Sciences of Ukraine, 3, Tereshenkivska st., 01004, Kyiv, Ukraine.

\vskip 3 mm

E-mail: \href{mailto:Yan.Sergiy@gmail.com}{Yan.Sergiy@gmail.com}

\end{document}